\documentclass{article}
\usepackage{amssymb}
\usepackage{graphicx}
\usepackage{amsmath}

\newtheorem{theorem}{Theorem}

\newtheorem{lemma}[theorem]{Lemma}

\newtheorem{problem}[theorem]{Problem}

\begin{document}

\title{Non-abelian free groups admit non-essentially free actions on rooted trees}
\author{Mikl\'{o}s Ab\'{e}rt and G\'{a}bor Elek}
\maketitle

\begin{abstract}
We show that every finitely generated non-abelian free group $\Gamma $
admits a spherically transitive action on a rooted tree $T$ such that the
action of $\Gamma $ on the boundary of $T$ is not essentially free. This
reproves a result of Bergeron and Gaboriau.

The existence of such an action answers a question of Grigorchuk,
Nekrashevich and Sushchanskii.
\end{abstract}

\section{Introduction}

Let $\Gamma $ be a finitely generated group. A \emph{chain} in $\Gamma $ is
a sequence $\Gamma =\Gamma _{0}\geq \Gamma _{1}\geq \ldots $ of subgroups of
finite index in $\Gamma $. Let $T=T(\Gamma ,(\Gamma _{i}))$ denote the \emph{%
coset tree}, a rooted tree on the set of right cosets of the subgroups $%
\Gamma _{n}$ defined by inclusion. The group $\Gamma $ acts by automorphisms
on $T$; this action extends to the boundary $\partial T$ of $T$, the set of
infinite rays starting from the root. The boundary is naturally endowed with
the product measure coming from the tree and $\Gamma $ acts by
measure-preserving homeomorphisms on $\partial T$. We call this action the 
\emph{boundary representation of }$\Gamma $ with respect to the chain $%
(\Gamma _{n})$. It is easy to see that this action is always ergodic and
minimal.

We say that the action of $\Gamma $ on $\partial T(\Gamma ,(\Gamma _{i}))$
is essentially free, if almost every element of $\partial T(\Gamma ,(\Gamma
_{i}))$ has trivial stabilizer in $\Gamma $. This is the case e.g. when the
chain consists of normal subgroups of $\Gamma $ and their intersection is
trivial.

The main aim of this note is to construct faithful non-essentially free
boundary representations of non-abelian free groups. Let $F_{d}$ denote the
free group of rank $d$.

\begin{theorem}
\label{main}Let $d\geq 2$ and let $\Gamma =F_{d}$. Then there exists a chain 
$(\Gamma _{n})$ in $\Gamma $ such that $\cap \Gamma _{n}=1$ and the boundary
representation of $\Gamma $ with respect to the chain $(\Gamma _{n})$ is not
essentially free.
\end{theorem}

This allows one to answer the question \cite[Problem 7.3.3]{grineksus}. We
say that two graphs are locally isomorphic if you can not distinguish them
using just local information, that is, the isomorphism classes of finite
balls coincide for the two graphs.

\begin{problem}[Grigorchuk, Nekrashevich and Sushchanskii]
Does there exist a spherically transitive group of automorphisms of a rooted
tree such that the Schreier graph of orbits, on the boundary, which are
typical in the sense of Baire category (the orbits of generic points), are
different from (not locally isomorphic to) the Schreier graph of orbits, on
the boundary, which are typical in the sense of measure?
\end{problem}

Indeed, it turns out that for our action of $F_{d}$, the Schreier graph of
an orbit of a Baire generic point is a $2d$-regular tree, while for a
measure generic point it is not.

\bigskip

\noindent \textbf{Remark. }The authors got aware that Theorem \ref{main} has
been proved by Bergeron and Gaboriau (see \cite[Theorem 4.1, point 5)]
{berggab}). Although the method of our proof is different, the result is the
same and therefore we shall leave this paper as an expository article.

\section{Proofs}

We start with some definitions. Let $(\Gamma _{n})$ be a chain in $\Gamma $.
Then the \emph{coset tree} $T=T(\Gamma ,(\Gamma _{n}))$ of $\Gamma $ with
respect to $(\Gamma _{n})$ is defined as follows. The vertex set of $T$
equals 
\begin{equation*}
T=\left\{ \Gamma _{n}g\mid n\geq 0,g\in \Gamma \right\}
\end{equation*}
and the edge set is defined by inclusion, that is, 
\begin{equation*}
(\Gamma _{n}g,\Gamma _{m}h)\text{ is an edge in }T\text{ if }m=n+1\text{ and 
}\Gamma _{n}g\supseteq \Gamma _{m}h\text{.}
\end{equation*}
Then $T$ is a tree rooted at the vertex $\Gamma $ and every vertex of level $%
n$ has the same number of children, equal to the index $\left| \Gamma
_{n}:\Gamma _{n+1}\right| $. The right actions of $\Gamma $ on the coset
spaces $\Gamma /\Gamma _{n}$ respect the tree structure and so $\Gamma $
acts on $T$ by automorphisms. This action is called the \emph{tree
representation} of $\Gamma $ with respect to $(\Gamma _{n})$.

The boundary $\partial T$ of $T$ is defined as the set of infinite rays
starting from the root. The boundary is naturally endowed with the product
topology and product measure coming from the tree. More precisely, for $%
t=\Gamma _{n}g\in T$ let us define the \emph{shadow} of $t$ as 
\begin{equation*}
\mathrm{Sh}(t)=\left\{ x\in \partial T\mid t\in x\right\}
\end{equation*}
the set of rays going through $t$. Set the base of topology on $\partial T$
to be the set of shadows and set the measure of a shadow to be 
\begin{equation*}
\mu (\mathrm{Sh}(t))=1/\left| \Gamma :\Gamma _{n}\right| .
\end{equation*}
This turns $\partial T$ into a totally disconnected compact space with a
Borel probability measure $\mu $. The group $\Gamma $ acts on $\partial T$
by measure-preserving homeomorphisms; we call this action the \emph{boundary
representation of }$\Gamma $ with respect to $(\Gamma _{n})$.

\bigskip

There are various levels of faithfulness of a boundary representation. Let 
\begin{equation*}
\partial T_{free}=\left\{ x\in \partial T\mid Stab_{\Gamma }(x)=1\right\} 
\text{.}
\end{equation*}
We say that the action is \emph{free}, if $\partial T_{free}=\partial T$.
The action is \emph{essentially free} (or that the chain satisfies the
Farber condition), if $\mu (\partial T\backslash \partial T_{free})=0$. The
action is \emph{topologically free} if $\partial T\backslash \partial
T_{free}$ is meager, i.e., a countable union of nowhere dense closed sets.
Note that the Farber condition has been introduced by Farber in \cite{farber}
in another equivalent formulation (the name `Farber condition' is from \cite
{berggab}).

It is easy to see that the following implications hold for a boundary
representation of a countable group $\Gamma $: 
\begin{eqnarray*}
(\Gamma _{n})\text{ is normal and }\left( \cap \Gamma _{n}=1\right)
&\Longrightarrow &\text{free}\Longrightarrow \text{essentially free}%
\Longrightarrow \\
&\Longrightarrow &\text{topologically free}\Longleftrightarrow \partial
T_{free}\neq \emptyset \Longrightarrow \text{faithful}
\end{eqnarray*}

For all but the third arrow it is easy to find examples showing that the
reverse implications do not hold.

\bigskip

Now we will start building towards Theorem \ref{main}. The first lemma is
straightforward from the definitions above.

\begin{lemma}
\label{fixratio}Let $\Gamma $ be a countable group with a chain $(\Gamma
_{n})$. Let 
\begin{equation*}
\mathrm{fixr}(g,\Gamma /\Gamma _{n})=\frac{\left| \left\{ x\in \Gamma
/\Gamma _{n}\mid xg=x\right\} \right| }{\left| \Gamma :\Gamma _{n}\right| }
\end{equation*}
denote the ratio of fixed points of $g$ acting on the right coset $\Gamma
/\Gamma _{n}$. Then the boundary representation of $\Gamma $ with respect to 
$(\Gamma _{n})$ is essentially free if and only if for all $g\in \Gamma $
with $g\neq 1$, the limit 
\begin{equation*}
\lim_{n\rightarrow \infty }\mathrm{fixr}(g,\Gamma /\Gamma _{n})=0\text{.}
\end{equation*}
\end{lemma}

\bigskip

Let $X$ be a set of symbols such that for all $x\in X$ the symbol $%
x^{-1}\notin X$. By an $X$\emph{-labeled graph} $G=(V,E,l)$ we mean a
finite, directed connected graph with vertex set $V$ and edge set $E$
together with a labeling function $l:E\rightarrow X$ that satisfies the
following: 
\begin{equation*}
\text{for all }v\in V,x\in X\text{ there is at most one }e\in E\text{
starting at }v\text{ such that }l(e)=x\text{.}
\end{equation*}

and 
\begin{equation*}
\text{for all }v\in V,x\in X\text{ there is at most one }e\in E\text{ ending
at }v\text{ such that }l(e)=x\text{.}
\end{equation*}
Note that we allow multiple edges and loops in $G$.

Let $G=(V,E,l)$ be an $X$-labeled graph. Then for every $x\in X$ we can
associate a function $f_{x}:V\rightarrow V$ as follows. Let $G^{\prime
}=(V,E^{\prime })$ be the directed graph obtained by erasing all edges from $%
G$ that is not $x$-labeled. For $v\in V$ let $C(v)$ denote the connected
component of $v$ in $G^{\prime }$. Then $C(v)$ is either an isolated point,
a directed circle or a directed simple path. If $v$ is an isolated point
then let $f_{x}(v)=v$. If $C(v)$ is a directed circle then let $%
f_{x}(v)=v^{\prime }$ where $(v,v^{\prime })\in E^{\prime }$. Similarly, if $%
C(v)$ is a path and $v$ is not the last point of $C(v)$ then let $%
f_{x}(v)=v^{\prime }$ where $(v,v^{\prime })\in E^{\prime }$. Finally, if $%
C(v)$ is a path and $v$ is the last point of $C(v)$ then let $%
f_{x}(v)=v^{\prime }$ where $v^{\prime }$ is the first point of $C(v)$.

It is easy to see that $f_{x}$ is a bijection of $V$ for all $x\in X$. Let $%
F_{X}$ denote the free group generated by the alphabet $X$. Then the mapping 
\begin{equation*}
\Phi :x\longmapsto f_{x}
\end{equation*}
extends to a homomorphism from $F_{X}$ to the symmetric group $\mathrm{Sym}%
(V)$, that is, a permutation action of $F_{X}$ on $V$. We will use the
following property of $\Phi $ -- the proof is a straightforward induction.

\begin{lemma}
\label{trivi}Let $w\in F_{X}$ be a reduced word of length $k$. Write 
\begin{equation*}
w=u_{1}u_{2}\cdots u_{k}
\end{equation*}
where $u_{i}\in X$ or $u_{i}^{-1}\in X$. Let $v_{0}\in V$ and for $i\geq 1$
let us define recursively $v_{i}\in V$ to be the vertex that satisfies 
\begin{equation*}
\begin{array}{cc}
l(v_{i-1},v_{i})=u_{i} & \text{if }u_{i}\in X \\ 
l(v_{i},v_{i-1})=u_{i}^{-1} & \text{if }u_{i}^{-1}\in X
\end{array}
\end{equation*}
Assume that the above recursive definition makes sense. Then 
\begin{equation*}
v_{k}=v_{0}\Phi (w)\text{.}
\end{equation*}
\end{lemma}

\bigskip

We are ready to prove Theorem \ref{main}.

\bigskip

\noindent \textbf{Proof of Theorem \ref{main}.}

Let 
\begin{equation*}
X=\left\{ a,b,c_{1},\ldots ,c_{d-2}\right\}
\end{equation*}
be an alphabet of $d$ letters. Let $\Gamma =F_{X}$ be the free group on $X$.
Let $C$ be the set of conjugacy classes in $\Gamma $ and let 
\begin{equation*}
C^{\prime }=\left\{ t\in C\mid \text{there is }w\in t\text{ starting with }%
a\right\}
\end{equation*}
Let us list the elements of $C^{\prime }$ as $t_{1},t_{2},\ldots $ and let $%
w_{i}\in t_{i}$ be a representative that starts with $a$ ($i\geq 1$). Let $%
k_{i}=\left| w_{i}\right| $ be the length of $w_{i}$ ($i\geq 1$) and let us
decompose 
\begin{equation*}
w_{i}=u_{i,1}u_{i,2}\cdots u_{i,k_{i}}
\end{equation*}
where $u_{i,j}\in X$ or $u_{i,j}^{-1}\in X$ ($1\leq j\leq k_{i}$). We can
assume that the sequence $k_{i}$ is non-decreasing. Let $0<\alpha <1\,$and
let $p_{i}$ be an increasing sequence of prime numbers satisfying 
\begin{equation*}
\prod_{i=1}^{n}(1-\frac{k_{i}+1}{p_{i}+k_{i}})>\alpha
\end{equation*}
for all $n\geq 1$.

\begin{figure}[tbp]
\includegraphics[width=12cm]{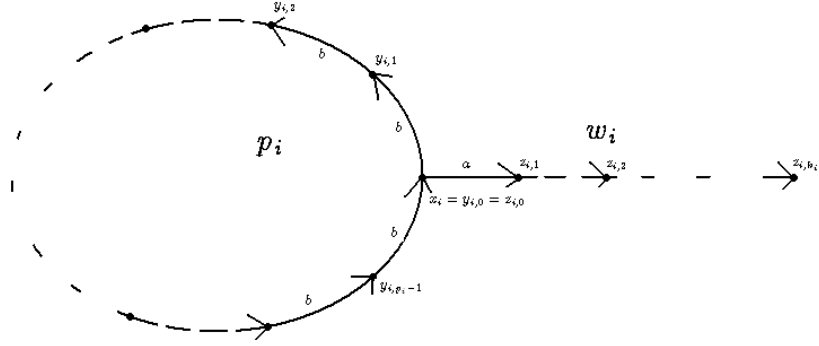}\newline
\caption{The graph $G_i$}
\label{fig:mining}
\end{figure}

For $i\geq 1$ let $G_{i}=(V_{i},E_{i},l_{i})$ be an $X$-labeled graph
defined as follows. Let 
\begin{equation*}
Y_{i}=\left\{ y_{i,0},y_{i,1},\ldots ,y_{i,p_{i}-1}\right\} \text{ and }%
Z_{i}=\left\{ z_{i,0},z_{i,1},\ldots ,z_{i,k_{i}}\right\}
\end{equation*}
be sets such that $x_{i}=y_{i,0}=z_{i,0}$ and $Y_{i}\cap Z_{i}=\left\{
x_{i}\right\} $. For convenience, denote $y_{i,p_{i}}=y_{i,0}$. Let $%
V_{i}=Y_{i}\cup Z_{i}$ and assume that $V_{i}\cap V_{j}=\emptyset $ ($i\neq
j $). Now for $0\leq j\leq p_{i}-1$ let 
\begin{equation*}
(y_{i,j},y_{i,j+1})\in E_{i}\text{ with }l_{i}(y_{i,j},y_{i,j+1})=b
\end{equation*}
and for $0\leq j\leq k_{i}-1$ let 
\begin{equation*}
\begin{array}{cc}
(z_{i,j},z_{i,j+1})\in E_{i}\text{ with }l_{i}(z_{i,j},z_{i,j+1})=u_{i,j+1}
& \text{if }u_{i,j+1}\in X \\ 
(z_{i,j+1},z_{i,j})\in E_{i}\text{ with }%
l_{i}(z_{i,j+1},z_{i,j})=u_{i,j+1}^{-1} & \text{if }u_{i,j+1}^{-1}\in X
\end{array}
\end{equation*}
Since $l_{i}(x_{i},y_{i,1})=b\neq a=l_{i}(x_{i},z_{i,1})$, and every other
vertex has in- and out-degree at most $1$, $G_{i}$ is indeed an $X$-labeled
graph.

Let $\Phi _{i}$ be the action of $F_{X}$ on $V_{i}$ defined by $G_{i}$. Then
the assumptions of Lemma \ref{trivi} hold for $G_{i}$ and $w_{i}$ with $%
v_{0}=x_{i}$, so we get 
\begin{equation*}
x_{i}\Phi _{i}(w_{i})=z_{i,k_{i}}\neq x_{i}
\end{equation*}

For $n\geq 1$ let 
\begin{equation*}
U_{n}=\bigoplus_{i=1}^{n}V_{i}\text{, }o_{n}=\bigoplus_{i=1}^{n}x_{i}\text{
and }\Psi _{n}=\bigoplus_{i=1}^{n}\Phi _{i}\text{ }
\end{equation*}
Let $O_{n}$ be the orbit of $o_{n}$ under $\Psi _{n}$ and let $H_{n}$ be the
stabilizer of $o_{n}$ in $F_{X}$. Then $(H_{n})$ is a chain in $F_{X}$ and
the right coset action of $F_{X}$ on $F_{X}/H_{n}$ is equal to the
restriction of $\Psi _{n}$ to $O_{n}$.

Let 
\begin{equation*}
U_{n}=\bigoplus_{i=1}^{n}Y_{i}\text{ and }P_{n}=\bigoplus_{i=1}^{n}\left(
Y_{i}\backslash \{x_{i}\}\right)
\end{equation*}
We claim that $U_{n}\subseteq O_{n}$. Indeed, $o_{n}\in U_{n}$ and $b\in
F_{X}$ acts on $Y_{i}$ as a cycle of length $p_{i}$. Since $\gcd
(p_{i},p_{j})=1$ ($i\neq j$), the orbit of $o_{n}$ under the cyclic group $%
\left\langle b\right\rangle $ equals $U_{n}$. By the definition of $\Phi
_{i} $, $a$ fixes $P_{n}$ pointwise, implying 
\begin{equation*}
\mathrm{fixr}(a,\Gamma /H_{n})\geq \frac{\left| P_{n}\right| }{\left|
O_{n}\right| }\geq \frac{\prod_{i=1}^{n}(p_{i}-1)}{%
\prod_{i=1}^{n}(p_{i}+k_{i})}=\prod_{i=1}^{n}(1-\frac{k_{i}+1}{p_{i}+k_{i}}%
)>\alpha
\end{equation*}
Using Lemma \ref{fixratio} this implies that the boundary representation of $%
\Gamma $ with respect to $(H_{n})$ is not essentially free.

Of course, $\cap _{n}H_{n}$ is not necessarily trivial. In fact, if $d>2$
then $c_{1}$ fixes $o_{n}$ ($n\geq 1$) so $c_{1}\in \cap _{n}H_{n}$. We
claim however, that the boundary representation of $\Gamma $ with respect to 
$(H_{n})$ is faithful. Assume it is not. Let $w\in \Gamma $ ($w\neq 1$) be
an element of the kernel. Then either $w$ or $w^{-1}$ is conjugate to $w_{m}$
for some $m$. But then $x_{m}\Phi _{m}(w_{m})\neq x_{m}$, implying $%
o_{m}\Psi _{m}(w_{m})\neq o_{m}$. So $w_{m}$ is not in the kernel of $\Psi
_{m}$, a contradiction. The claim holds.

Now we invoke a result in \cite{chains} saying that every faithful boundary
representation of a countable free group is topologically free. We get that
there exists $x\in \partial T$ with trivial stabilizer in $\Gamma $. Let $%
t_{i}$ be the vertex lying on the ray $x$ of level $i$ and let $\Gamma _{i}=%
\mathrm{Stab}_{\Gamma }(t_{i})$. Then $\cap _{i}\Gamma _{i}=1$ and since $%
\Gamma $ acts transitively on each level of $T$, $\Gamma _{i}$ is conjugate
to $H_{i}$ in $\Gamma $ ($i\geq 0$). This implies that the boundary action
of $\Gamma $ with respect to $(\Gamma _{i})$ is isomorphic to the boundary
action with respect to $(H_{i})$. In particular, the boundary action of $%
\Gamma $ with respect to $(\Gamma _{i})$ is not essentially free. The
theorem is proved. $\square $

\bigskip

\noindent \textbf{On a problem of Grigorchuk, Nekrashevich and Sushchanskii.}
Now we show how this leads to the solution of \cite[Problem 7.3.3]{grineksus}%
. Fix $d\geq 2$ and a minimal generating set $X$ in $F_{d}$. Let us take the
chain constructed in Theorem \ref{main}\textbf{. }Then the action of $F_{d}$
on the coset tree is spherically transitive and topologically free, that is,
a Baire generic point of the boundary has trivial stabilizer in $F_{d}$.
This implies that the Schreier graph of an orbit of a Baire typical point is
isomorphic to the Cayley graph $\mathrm{Cay}(F_{d},X)$, which is an infinite 
$2d$-regular tree. On the other hand, the boundary representation is not
essentially free, which implies that the Schreier graph of an orbit of a
measure typical point is isomorphic to the Schreier graph $\mathrm{Sch}%
(F_{d}/H,X)$ where $H$ is a nontrival subgroup of $F_{d}$. This graph is
hence also $2d$-regular but is never a tree. Thus there exists a ball in $%
\mathrm{Sch}(F_{d}/H,X)$ which can not be embedded into $\mathrm{Cay}%
(F_{d},X)$.

\bigskip

\noindent \textbf{Remark. }Our method provides a chain in $F_{n}$ that has
rapidly growing index or, equivalently, the number of children of a vertex
on the coset tree grows very fast. However, there is no reason to assume
that this is necessary. In fact, an infinite rooted binary tree should be
sufficient. It may be possible to use a random method to show that an
odometer (a.k.a. adding machine) and a random element fixing a fixed nowhere
dense set of positive measure on the boundary of a rooted binary tree
generates a free group a.s.

\end{document}